\documentclass[10pt]{amsart}
\usepackage{latexsym}
\usepackage{amsmath}
\usepackage{amssymb}
\usepackage{mathrsfs}
\usepackage{graphicx}
\usepackage{xcolor}
\usepackage{bbm}
\usepackage{bbold}
\input diagxy.tex
\def\epi{\mathop{\fam0 epi}\nolimits}
\def\cop{\mathop{\fam0 cop}\nolimits}
\def\dom{\mathop{\fam0 dom}\nolimits}
\def\On{\mathop{\fam0 On}\nolimits}
\def\scyc{\mathop{\fam0 mix}\nolimits}
\def\st{\operatorname{st}}
\def\ltd{\operatorname{ltd}}

\def\bnorml{\mathopen{\kern1pt \vrule height6.5pt depth1.5pt width1pt\kern1.5pt}}
\def\bnormr{\mathclose{\kern1pt \vrule height6.5pt depth1.5pt width1pt\kern1pt}}

\let\subsection\relax%\newpage
\let\cite\relax
\let\footnote\relax

\thanks
{This article is an~extended version of a~talk prepared for
the International
Conference ``Methods of Logic in Mathematics V,''
June 1--6, 2008, St. Petersburg}

\begin{document}
\title{Nonstandard Models and Optimization}
\author{S. S. Kutateladze}
\date{May 26, 2008}
\begin{abstract}
This is an overview of a few possibilities that are open by
model theory in applied mathematics.  Most attention is paid to
the present state and frontiers of the Cauchy method of majorants,
approximation of operator equations with finite-dimensional analogs, and the Lagrange multiplier principle
in multiobjective decision making.
\end{abstract}
\address[]{
 Sobolev Institute of Mathematics\newline
%\indent Universitetski\u\i{} pr.~4\newline
\indent Novosibirsk
}
\email{
sskut@math.nsc.ru
}

\maketitle
\section{\large\bf\color{blue} Agenda}
The union of functional analysis and applied mathematics
celebrates its sixtieth anniversary
this year. %\footnote{Cp.~\cite{KantLGU, KantUMN}.}
This talk focuses on the trends of interaction between
model theory and the methods of domination, discretization, and scalarization.

\section{\large\bf\color{blue}The  Art of Calculus}
Provable counting is the art of calculus which is mathematics
in modern parlance. Mathematics exists as a science
more than two and a half millennia, and we can never mixed
it with history or chemistry. In this respect
our views of what is mathematics are independent of time.

The objects of mathematics  are the quantitative
forms of human reasoning. Mathematics functions as the science
of convincing calculations. Once-demonstrat\-ed, the facts of mathematics
well never vanish. Of course, mathematics renews itself constantly,
while the stock increases of mathematical notions and construction and
the understanding changes of the rigor and technologies of proof and
demonstration. The frontier we draw between  pure
and applied  mathematics is also time-dependent.

\section{\large\bf\color{blue}Francis Bacon}

{\color{red}\it The Mathematics are either pure or mixed. To the Pure
Mathematics are those sciences belonging which handle quantity
determinate, merely severed from any axioms of natural philosophy; and
these are two, Geometry and Ari\-thmetic; the one handling quantity
continued, and the other dissevered. Mixed hath for subject some
axioms or parts of natural philosophy, and considereth quantity
determined, as it is auxiliary and incident unto them\dots.
%For many parts
%of nature can neither be invented with sufficient subtilty, nor
%demonstrated with sufficient perspicuity, nor accommodated unto use
%with sufficient dexterity, without the aid and intervening of the
%mathematics; of which sort are perspective, music, astronomy,
%cosmography, architecture, enginery, and divers others.

In the Mathematics\dots
% I can report no deficience, except it be that
%men do not sufficiently understand the excellent use of the Pure
%Mathematics, in that they do remedy and cure many defects in the wit
%and faculties intellectual\dots.
%For if the wit be too dull, they sharpen
%it; if too wandering, they fix it; if too inherent in the sense, they
%abstract it. So that as tennis is a game of no use in itself, but of
%great use in respect it maketh a quick eye and a body ready to put
%itself into all postures; so
%in the Mathematics,
that use which is
collateral and intervenient is no less worthy than that which is
principal and intended\dots. And as for the Mixed Mathematics, I may
only make this prediction, that there cannot fail to be more kinds of
them, as nature grows further disclosed.}

{\small\hfill The Advancement of Learning, 1605}
%\footnote{The complete
%title was as follows: ``The tvvoo bookes of Francis Bacon, of the
%proficience and aduancement of learning, diuine and humane. To the
%King. At London: Printed for Henrie Tomes, 1605.''}

\section{\large\bf\color{blue}Mixed Turns into Applied}
After the lapse of 150 years   Leonhard Euler used the words
``pure  mathematics'' in the title of one of his papers
{\it Specimen de usu
observationum in mathesi pura\/} in~1761. It was practically at the same time
that the term ``pure mathematics'' had appeared in the
eldest  {\it Encyclopaedia Britannica}. In the nineteenth century
``mixed'' mathematics became to be referred to as ``applied.''

The famous  {\it Journal de Math\'ematiques Pures et Appliqu\'ees}
was founded by Joseph Liouville in~1836 and {\it The
Quarterly Journal of Pure and Applied Mathematics\/} started publication
in 1857.

\section{\large\bf\color{blue}Pure and Applied Mathematics}
The intellectual challenge, beauty, and intrinsic
logic of the topics under study are the impetus of many comprehensive and
deep studies in  mathematics which are customarily qualified as pure.
%Knowledge of the available mathematical methods and the understanding
%of their power underlie the applications of mathematics in other sciences.
Any application of mathematics is impossible without creating some
metaphors, models of the phenomena and processes under examination.
Modeling is a special independent sphere
of intellectual activities which is out of mathematics.

Application of mathematics resides beyond ma\-thematics
in much the same way as maladies exist in nature rather than inside
medicine. Applied mathematics acts as an apothecary
mixing drugs for battling illnesses.

The art and craft of mathematical techniques
for the problems of other sciences are the content of
applied mathematics.

\section{\large\bf\color{blue}New Challenges}

Classical mechanics in the broadest sense of the words was the
traditional sphere of applications of  mathematics in the
nineteenth century.The beginning of the twentieth century was marked with a sharp
enlargement of the sphere of applications of mathematics.
Quantum mechanics appeared, requesting for new mathematical tools.
The theory of operators in Hilbert spaces and distribution theory
were oriented   to adapting the heuristic methods
of the new physics. At the same time the social phenomena became the
object of the nonverbal research requiring the invention of
especial mathematical methods.
The demand for the statistical treatment of various data grew rapidly.
Founding new industries as well as introducing of promising
technologies and new materials, brought about the necessity of
elaboration of the technique of calculations.
The rapid progress of applied mathematics was facilitated
by the automation and mechanization of accounting and standard calculations.

\section{\large\bf\color{blue}Cofathers of New Mentality}
In the 1930s applied mathematics  rapidly approached
functional analysis.

Of profound importance in this trend
was the research of John von Neumann in the mathematical foundations
of quantum mechanics and game theory as a tool for economic
studies.

Leonid Kantorovich was a pioneer and generator of new synthetic ideas
in Russia.

\section{\large\bf\color{blue}Enigmas of Economics}
The main particularity of the extremal problems of economics
consists in the presence of numerous conflicting ends and interests
to be harmonized. We  encounter the instances of multicriteria  optimization. Seeking for an optimal solution in these circumstances, we
must take into account various contradictory preferences which combine
into  a sole compound aim.  It is impossible as a rule
to distinguish some particular  scalar target and ignore the rest of
the targets. This circumstance involves the specific difficulties
that are untypical in the scalar case: we must specify
 what we should call a solution of a vector program and we must agree
 upon the method of conforming versatile ends provided that some
 agreement is possible  in principle.
Therefore, it is actual to seek for the reasonable
concepts of optimality in multiobjective decision making.
Among these we distinguish  the concepts of ideal and generalized
optimum alongside  Pareto-optimum as well as approximate and
infinitesimal optimum.

\section{\large\bf\color{blue}Enter the Reals}
Optimization is the science of choosing the best.
To choose, we use preferences. To optimize, we use infima and
suprema (for bounded subsets) which is practically
the {least upper bound property}.
So optimization needs ordered sets and primarily
(boundedly) complete lattices.

To operate with preferences,
we use group structure. To aggregate and scale, we use linear structure.

All these are happily provided by the {\it reals\/} $\mathbb R$, a one-dimensional
Dedekind complete vector lattice. A Dedekind  complete vector
lattice is a {\it Kantorovich space}.

\section{\large\bf\color{blue}Scalarization}

Scalarization in the most general sense means reduction to numbers.
Since each number is a measure of quantity,
the idea of scalarization is clearly of a universal importance
to mathematics. The deep roots of scalarization are revealed by
the Boolean valued validation of the Kantorovich heuristic principle.
We will dwell upon the aspects of scalarization most important in applications
and connected with the problems of multicriteria optimization.

\section{\large\bf\color{blue}Legendre in Disguise}
Assume that $X$~is a vector space, $E$ is an~ordered vector space,
$f:X\rightarrow E^\bullet:=E\cup{+\infty}$ is a~convex operator,
and $C:=\dom(f)\subset X$ is a~convex set.
 A {\it vector program\/}
 $(C,f)$  is written  as follows:
  $$
  x\in C,\ \ f(x)\rightarrow \inf\!.
 $$

The standard sociological trick  includes $(C,f)$
into a~parametric family  yielding
the  {\it Legendre trasform\/} or {\it Young--Fenchel transform\/} of $f$:
$$
f^*(l):=\sup_{x\in X}{(l(x)-f(x))},
$$
with $l\in X^{\#}$  a linear functional over $X$.
The epigraph of $f^*$  is a convex subset of
$ X^{\#}$ and so $f^*$  is  convex.
Observe that $-f^*(0)$ is the value of $(C,f)$.

\section{\large\bf\color{blue}Order Omnipresent}
A convex function is locally a positively homogeneous convex
function, a {\it sublinear functional}. Recall that
$p: X\to\mathbb R$ is sublinear whenever
$$
\epi p:=
\{ (x,\ t)\in X\times\mathbb R: p(x)\le t\}
$$
is a cone. Recall that a numeric function is uniquely
determined from its epigraph.

Given
$C\subset X$, put
$$
H(C):=\{(x,\ t)\in X\times\mathbb R^+ : x\in tC\},
$$
the  {\it H\"ormander transform\/} of $C$.
Now, $C$ is convex if and only if $H(C)$  is a~cone.
A~space with a cone is a {\it $($pre$)$ordered vector space}.

{\color{red}\it The order, the symmetry, the harmony enchant us\dots.}

{\small\hfill Leibniz}

\section{\large\bf\color{blue}Fermat's Criterion}
$\partial f(\bar x)$, the {\it subdifferential\/} of $f$
at  $\bar x$, is
$$
\{l\in X^\# :(\forall x\in X)\
l(x)-l(\bar x)\leq f(x)-f(\bar x) \}.
$$

A point $\bar x$ is a solution
to the  minimization problem $(X,f)$ if and only if
$$
0\in \partial f(\bar x).
$$

This {\it Fermat criterion} turns into the Rolle Theorem in a smooth case
and is  of little avail without effective tools for calculating $\partial f(\bar x).$
A~convex analog of the  ``chain rule''  is in order.

\section{\large\bf\color{blue}Enter Hahn--Banach}
The {\it  Dominated Extension} takes the form
$$
\partial (p\circ \iota) (0)=(\partial p)(0)\circ \iota,
$$
with $p$ a~sublinear functional over $X$ and
$\iota $ the identical embedding of some~subspace of~$X$ into~$X$.

If the target $\mathbb R$ may be replaced with an ordered vector space
$E$, then $E$ admits {\it dominated extension}.

%In fact, the situation is even more spectacular:

%{\bf Theorem.}
%{\sl A preordered vector space $E$ admits dominated
%extension of linear operators if and only if $E$ has the least upper
%bound property}.

\section{\large\bf\color{blue}Enter Kantorovich}

The matching of convexity and order
was established in two steps.

{\bf Hahn--Banach--Kantorovich  Theorem.}
{\sl Every Kantorovich space admits dominated extension of linear operators}.

This theorem proven by Kantorovich in~1935 was a first attractive result
of the theory of ordered vector spaces.

{\bf Bonnice--Silvermann--To Theorem.}
{\sl Each ordered vector space
admitting dominated extension of linear operators is a~Kantorovich space}.

\section{\large\bf\color{blue}New Heuristics}
Kantorovich demonstrated the role of  $K$-spaces by the example
of the Hahn--Banach theorem. He proved that this central principle of functional analysis admits the replacement of reals with elements of an arbitrary  $K$-space while substituting linear and sublinear operators with range
in this space for linear and sublinear functionals.
These observations laid grounds for the universal heuristics based on
his intuitive belief that the members of an abstract Kantorovich space
are a sort of generalized numbers.

\section{\large\bf\color{blue}Canonical Operator}

Consider a~Kantorovich space
$E$
and an arbitrary nonempty set
$\mathfrak A$.
Denote by
$l_\infty (\mathfrak A,E)$
the set of all order bounded mappings from
$\mathfrak A$
into $E$; i.e.,
$f\in l_\infty (\mathfrak A,E)$
if and only if
$f:\mathfrak A \to E$
and
$\{f(\alpha):\alpha\in\mathfrak A\}$
is order bounded in $E$.
It is easy to verify that
$l_\infty (\mathfrak A,E)$ becomes a Kantorovich
space if endowed with the coordinatewise algebraic operations and order.
The operator
$\varepsilon_{\mathfrak A, E}$
acting from
$l_\infty (\mathfrak A,E)$
into
$E$
by the rule
$$
\varepsilon_{\mathfrak A, E}:f\mapsto\sup \{f(\alpha):\alpha\in\mathfrak A\}
\quad (f\in l_\infty (\mathfrak A,E))
$$
is called the {\it canonical sublinear operator\/}
given
$\mathfrak A$
and
$E$.
We often write
$\varepsilon_{\mathfrak A}$
instead of
$\varepsilon_{\mathfrak A, E}$
when it is clear from the context what Kantorovich space is meant.
The notation
$\varepsilon_n$
is used when the cardinality of~$\mathfrak A$
equals
$n$ and we call the operator
$\varepsilon_n$
{\it finitely-generated}.

\section{\large\bf\color{blue}Support Hull}

Consider a~set~
$\mathfrak A$
of linear operators acting from
a~vector space
$X$
into a~Kantorovich space
$E$.
The set
$\mathfrak A$
is {\it weakly order  bounded\/} if
$\{\alpha x:\alpha\in\mathfrak A\}$
is order bounded for every
$x\in X$.
We denote by
$\langle\mathfrak A\rangle x$
the mapping that assigns the element
$\alpha x\in E$
to each
$\alpha\in \mathfrak A$,
i.e.
$\langle\mathfrak A\rangle x: \alpha\mapsto\alpha x$.
If
$\mathfrak A$
is weakly order bounded then
$\langle\mathfrak A\rangle x\in l_\infty (\mathfrak A,E)$
for every fixed
$x\in X$.
Consequently, we obtain the linear operator
$\langle\mathfrak A\rangle:X\rightarrow l_\infty (\mathfrak A,E)$
that acts as
$\langle\mathfrak A\rangle:x\mapsto\langle\mathfrak A\rangle x$.
Associate with
$\mathfrak A$
one more operator
$$
p_{\mathfrak A}: x\mapsto\sup \{\alpha x: \alpha\in\mathfrak A\}\quad(x\in X).
$$
The operator
$p_{\mathfrak A}$
is sublinear. The support set
$\partial p_{\mathfrak A}$
is denoted by
$\cop (\mathfrak A)$
and referred to as  the {\it support hull\/} of
$\mathfrak A$.

\section{\large\bf\color{blue}  Hahn--Banach in Disguise  }

{\bf Theorem.} {\sl If
$p$
is a~sublinear operator with
$\partial p=\cop (\mathfrak A)$
then $
P=\varepsilon_{\mathfrak A}\circ \langle\mathfrak A\rangle.
$
Assume further that
$p_1: X\to E$
is a~sublinear operator and
$p_2: E\to F$
is an increasing sublinear operator. Then
$$
\partial (p_2\circ p_1)
=%\left
\{ T\circ\langle\partial p_1\rangle:
 T\in L^+(l_{\infty}(\partial p_1, E),F)
\&\ T\circ\Delta_{\partial p_1}\in \partial p_2
%\right
\}.
$$
Moreover, if
$\partial p_1=\cop (\mathfrak A_1)$
and
$\partial p_2=\cop (\mathfrak A_2)$
then}
$$
\partial (p_2\circ p_1)
=\bigl\{T\circ\langle\mathfrak A_1\rangle : T\in L^+
(l_{\infty}(\mathfrak A_1,E),F)\
\left(\exists\alpha\in\partial\varepsilon_{\mathfrak A_2}\bigr)\
T\circ\Delta_{\mathfrak A_1}=\alpha\circ\langle\mathfrak A_2\rangle\right\}.
$$

%5Hahn--Banach in the classical formulation is of course the simplest
%chain rule for removing any linear embedding from the  subdifferential sign.

\section{\large\bf\color{blue}Enter Boole}
Cohen's final solution of the problem of the cardinality of the
continuum within ZFC gave rise to the Boolean-valued models
by Vop\v enka, Scott, and Solovay.
Takeuti coined the term ``Boolean-valued
analysis''  for applications of the new models
to functional analysis.

Let
$B$
be a~complete Boolean algebra. Given an ordinal
$\alpha$,
put
$$
V_{\alpha}^{(B)}
:=\{x:
(\exists \beta\in\alpha)\ x:\dom (x)\rightarrow
B\ \&\ \dom (x)\subset V_{\beta}^{(B)}  \}.
$$
The {\it Boolean-valued
universe\/}
${\mathbb V}^{(B)}$
is
$$
{\mathbb V}^{(B)}:=\bigcup\limits_{\alpha\in\On} V_{\alpha}^{(B)},
$$
with $\On$ the class of all ordinals. The truth
value $[\![\varphi]\!]\in B$ is assigned to each formula
$\varphi$ of ZFC relativized to ${\mathbb V}^{(B)}$.

\section{\large\bf\color{blue}~Enter Descent}

Given $\varphi$, a~formula of ZFC, and
$y$, a~subset ${\mathbb V}^{B}$; put
$A_{\varphi}:=A_{\varphi(\cdot,\ y)}:=\{x:\varphi (x,\ y)\}$.
The {\it descent\/}
$A_{\varphi}\!\!\downarrow$
of a~class
$A_{\varphi}$
is
$$
A_{\varphi}\!\!\downarrow:=\{t:t\in {\mathbb V}^{(B)} \ \&\  [\![\varphi  (t,\ y)]\!]=\mathbb 1\}.
$$
If
$t\in A_{\varphi}\!\!\downarrow$,
then
it is said
that
{\it $t$
satisfies
$\varphi (\cdot,\ y)$
inside
${\mathbb V}^{(B)}$}.

%The descent of each class is {\it strongly cyclic}, i.e. it  contains all mixings of its elements. Moreover,
%two classes inside
%${\mathbb V}^{(B)}$
%coincide if and only if they consist of the same elements inside
%${\mathbb V}^{(B)}$.

The {\it descent\/}
$x\!\!\downarrow$
of an element
$x\in {\mathbb V}^{(B)}$
is defined by the rule
$$
x\!\!\downarrow:=\{t: t\in {\mathbb V}^{(B)}\ \&\  [\![t\in x]\!]=\mathbb 1\},
$$
i.e. $x\!\!\downarrow=A_{\cdot\in x}\!\!\downarrow$.
The class $x\!\!\downarrow$ is a~set. Moreover,
$x\!\!\downarrow\subset\scyc (\dom (x))$,
where $\scyc$ is the symbol of the taking of the
{\it strong cyclic hull}.
If $x$ is a~nonempty set
inside
${\mathbb V}^{(B)}$ then
$$
(\exists z\in x\!\!\downarrow)[\![(\exists z\in x)\ \varphi (z)]\!]
=[\![\varphi(z)]\!].
$$

\section{\large\bf\color{blue} The Reals in Disguise}

There is an object
$\mathscr R$
inside
${\mathbb V}^{(B)}$ modeling $\mathbb R$, i.~e.,
$$
[\![\mathscr R\ {\text{is the reals}}\,]\!]=\mathbb 1.
$$
Let $\mathscr R\!\!\downarrow$ stand for the descent of
 the carrier $|\mathscr R|$ of the algebraic system
$\mathscr R:=(|\mathscr R|,+,\,\cdot\,,0,1,\le)$
inside ${\mathbb V}^{(B)}$.
Implement the descent of the structures on $|\mathscr R|$
to $\mathscr R\!\!\downarrow$ as follows:
$$
\gathered
x+y=z\leftrightarrow [\![x+y=z]\!]=\mathbb 1;
\\
xy=z\leftrightarrow [\![xy=z]\!]=\mathbb 1;
\\
x\le y\leftrightarrow [\![x\le y]\!]=\mathbb 1;
\\
\lambda x=y\leftrightarrow [\![\lambda^\wedge x=y]\!]=\mathbb 1
\\
(x,y,z \in\mathscr R\!\!\downarrow,\ \lambda\in\mathbb R).
\endgathered
$$

{\bf Gordon Theorem.}
{\sl $\mathscr R\!\!\downarrow$
with the descended structures is a universally complete
Kantorovich space with base
$\mathscr B(\mathscr R\!\!\downarrow)$
isomorphic to
$B$}.

\section{\large\bf\color{blue}Norming Sequences}

$$
\bnorml(\xi_1,\xi_2,\dots)\bnormr=
(|\xi_1|,|\xi_2|,\dots,|\xi_{N-1}|, \sup\limits_{k\geq N}|\xi_k|)\in\mathbb R^N.
$$

\begin{figure}[tbp]
\includegraphics[scale=1.0]{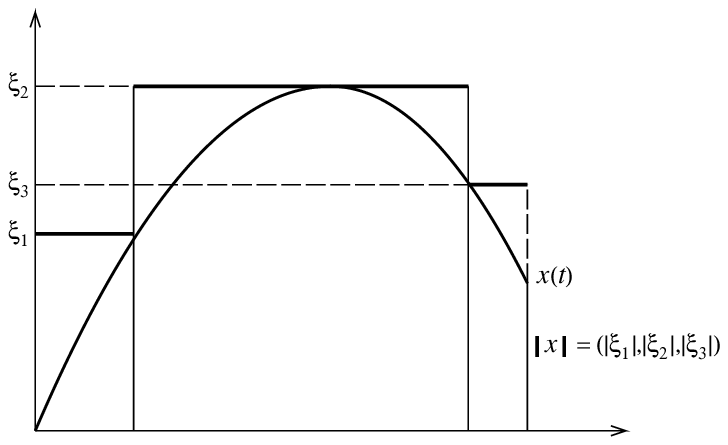}
\end{figure}

{\color{red}\it I believe that the use of
members of semi-ordered linear spaces instead of reals in various
estimations can lead to essential improvement of the latter.}

{\small\hfill Kantorovich, Herald of LGU, {\bf 6}, 3--18 (1948)}
%\centerline{\epsfxsize8cm\epsfbox{ris1.eps}}

\section{\large\bf\color{blue}Domination}

Let $X$ and $Y$ be real vector spaces lattice-normed
with $K$-spaces $E$ and $F$. In other words,
given are some lattice-norms ${\bnorml\cdot\bnormr}_{X}$
and ${\bnorml\cdot\bnormr}_{Y}$. Assume further that $T$
is a linear operator from~$X$ to~$Y$ and $S$ is a~positive
operator from  $X$ into~$Y$ satisfying
 $$
  \bfig
\square[X`Y`E`F;T`{\bnorml\cdot\bnormr}_{X}`{\bnorml\cdot\bnormr}_{Y}`S]
 \efig
 $$

Moreover, in case
$$
  {\bnorml Tx \bnormr}_{Y} \leq S{\bnorml x\bnormr}_{X}\quad (x\in X),
$$

\noindent
we call $S$  the~{\it dominant\/} or {\it majorant\/} of~$T$.

\section{\large\bf\color{blue}Enter Abstract Norm}
If the set of all dominants of~ $T$ has the least element, then the latter
is called the {\it abstract norm\/}  or {\it least dominant\/} of $T$ and denoted by
 $\bnorml T\bnormr $.
Hence, the least dominant  $\bnorml T\bnormr $~is the least
positive operator from~$E$ to~ $F$ such that
 $$
   \bnorml Tx\bnormr \leq \bnorml T\bnormr (\bnorml  x \bnormr )\quad (x\in X).
 $$

\section{\large\bf\color{blue}Domination and Model Theory}
These days the development of domination proceeds within the frameworks
of Boolean valued analysis.
All principal properties of lattice normed spaces
represents the Boolean valued interpretations of the relevant properties
of classical normed spaces. The most important interrelations here are
as follows: Each Banach space inside a Boolean valued model becomes
a universally complete Banach--Kan\-to\-rovich spaces in result of
the external deciphering of constituents. Moreover,
each lattice normed space  may be realized as a dense subspace
of some Banach space in an appropriate Boolean valued model.
Finally, a Banach space~$X$ results from some
Banach space inside  a Boolean valued model
by a special machinery of bounded descent
if and only if $X$ admits a complete Boolean
algebra of norm-one projections which enjoys the
cyclicity property. The latter amounts to the
fact that $X$ is a~Banach--Kantorovich space and  $X$
is furnished with a~mixed norm.

\section{\large\bf\color{blue}Approximation}

Convexity is an abstraction of finitely many stakes encircled with a surrounding rope, and so no variation of stakes can ever spoil the convexity
of the tract to be surveyed.

Study of stability in optimization
is accomplished sometimes by introducing various  epsilons
in appropriate places. One of the earliest excursions in this direction
is connected with the classical  Hyers--Ulam stability theorem for
$\varepsilon$-convex functions.
Exact calculations with epsilons and sharp estimates
are sometimes bulky and slightly mysterious.  Some alternatives
are suggested
by actual infinities, which is illustrated with the conception
of {\it infinitesimal optimality}.

\section{\large\bf\color{blue}Enter Epsilon and Monad}

Assume given a ~convex operator
$f:X\to E\cup{+\infty}$
and a~ point
$\overline x$
in the effective domain
$\dom(f):=\{x\in X:f(x)<+\infty\}$
of
~$f$.
Given
$\varepsilon \ge 0$
in the positive cone
$E_+$
of
$E$,
by the
$\varepsilon $-{\it subdifferential\/}
of~$f$
at
~$\overline x$
we mean the set
$$
\partial\, {}^\varepsilon\!f(\overline x)
:=\big\{T\in L(X,E):
(\forall x\in X)(Tx-Fx\le T\overline x -f\overline x+\varepsilon) \big\},
$$
with
$L(X,E)$
standing as usual for the space of linear operators
from~
$X$
to
~$E$.

Distinguish
some downward-filtered subset
~$\mathscr E$ of
$E$
that is composed of positive elements.
Assuming
$E$ and~$\mathscr E$
 standard, define the {\it monad\/}
$\mu (\mathscr E)$ of $\mathscr E$ as
$\mu (\mathscr E):=\bigcap\{[0,\varepsilon ]:\varepsilon \in
{}^\circ\!\mathscr E\}$.
The members of $\mu(\mathscr E)$ are {\it positive
infinitesimals\/}  with respect to~$\mathscr E $.
As usual,
${}^\circ\!\mathscr E$
denotes the external set of~
all standard members of
~$E$,
the {\it standard part\/} of
~$\mathscr E$.

\section{\large\bf\color{blue}Pareto Optimality}

\subsection{}
 Fix a positive element~$\varepsilon\in E$.
 A feasible point~$x_0$ is a
 {\it $\varepsilon$-solution\/}
or {\it $\varepsilon$-optimum\/}
of a program~$(C,f)$ provided that $f(x_0)\leq e+\varepsilon$ with
 $e$~the value of~$(C,f)$. In other words, $x_0$
 is an $\varepsilon$-solution of~$(C,f)$
if and only if $x_0\in C$ and the $f(x_0)-\varepsilon$~is the greatest
lower bound of~$f(C)$ or, equivalently,
 $f(C)+\varepsilon\subset f(x_0)+E^+$.
 Clearly,~$x_0$ is a $\varepsilon$-solution
of an unconditional problem~$f(x)\rightarrow \inf$ if and only if
the zero belong to~$\partial^\varepsilon f(x_0)$; i.~e.,
 $$
 f(x_0)\leq \inf\limits_{x\in X}\,f(x)+\varepsilon\,
 \leftrightarrow\, 0\in \partial_\varepsilon f(x_0).
 $$

\section{\large\bf\color{blue}Approximate Efficiency}
 A feasible point~$x_0$ is
 {\it $\varepsilon$-Pareto optimal} for $(C,f)$
whenever $f(x_0)$~is a minimal element of~$U+\varepsilon$,
with $U:=f(C)$; i.~e.,
 $(f(x_0)-E^+)\cap(f(C)+\varepsilon)=[f(x_0)]$.
In more detail,   $x_0$ is $\varepsilon$-Pareto-optimal
means that $x_0\in C$ and, for all $x\in C$,
from $f(x_0)\geq f(x)+\varepsilon$ it follows that
 $f(x_0)\sim f(x)+\varepsilon$.

\medskip
\begin{figure}
\includegraphics[scale=1.0]{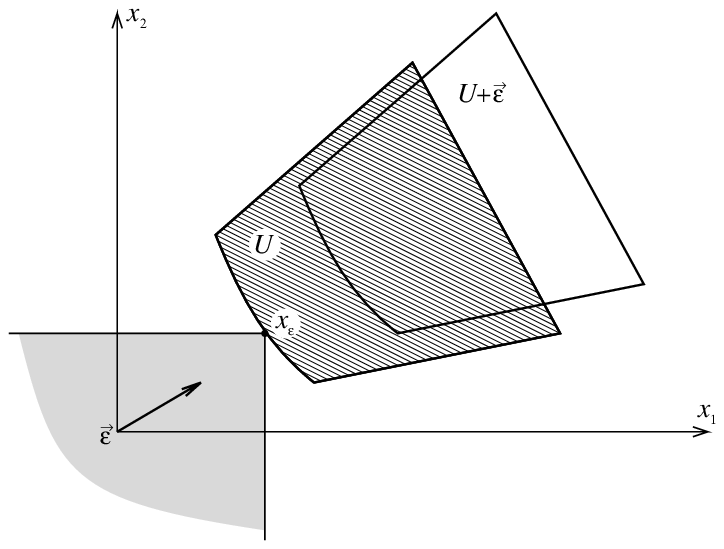}
\end{figure}
%\centerline{\epsfxsize8cm\epsfbox{ris2.eps}}

\section{\large\bf\color{blue}Subdifferential Halo}
Assume that the monad $\mu (\mathscr E )$
is an external cone over ${}^\circ \mathbb R $ and, moreover,
$\mu (\mathscr E)\cap{}^\circ\! E=0$.
In application, $\mathscr E $ is usually
the filter of order-units of $E$.
The relation of
{\it infinite proximity\/} or
{\it infinite closeness\/}
between the members of $E$ is introduced as follows:
$$
e_1 \approx e_2 \leftrightarrow e_1 -e_2 \in\mu
(\mathscr E )\ \&\ e_2 -e_1 \in\mu (\mathscr E ).
$$

Now
$$
Df(\overline x):=\bigcap\limits_{\varepsilon \in{}^\circ \mathscr E }\,
\partial ^\varepsilon f(\overline x)=
\bigcup\limits_{\varepsilon \in\mu (\mathscr E )}\,
\partial^\varepsilon f(\overline x);
$$
the {\it infinitesimal subdifferential} of
$f$ at  $\overline x$.
The elements of
$Df(\overline x)$ are
{\it infinitesimal subgradients\/}
of $f$ at
~$\overline x$.

\section{\large\bf\color{blue}Exeunt Epsilon}

{\scshape Theorem.} {\sl
Let $f_1:X\times Y\rightarrow E\cup +\infty$ and
$f_2:Y\times Z\rightarrow E\cup +\infty$ be convex operators.
Suppose that the convolution
$f_2\vartriangle f_1$ is infinitesimally exact at some point $(x,y,z)$; i.e.,
$
(f_2\vartriangle f_1)(x,y)\approx f_1(x,y)+f_2(y,z).
$
If, moreover, the
convex sets $\epi(f_1,Z)$ and $\epi(X,f_2)$ are
in general
position then}
$$
D(f_2\vartriangle f_1)(x,y)=
Df_2(y,z)\circ Df_1(x,y).
$$

\section{\large\bf\color{blue}Discretization}

{\color{red}\it It seems to me that the main idea of this theory is of a general character
and reflects the general gnoseological principle for studying
complex systems. It was, of course, used earlier, and it is also used
in systems analysis, but it does not have a rigorous mathematical apparatus.

The principle consists simply  in  the fact that
to a given large complex system in some space a simpler, smaller dimensional
model in this or a~simpler space is associated by means of one-to-one or
one-to-many  correspondence. The study of this simplified model
turns out, naturally, to be simpler and more practicable.
This method, of course, presents definite requirements on the quality
of the approximating system.}

{\small\hfill Kantorovich, Herald of LGU, {\bf 6}, 3--18 (1948)}

\section{\large\bf\color{blue}Hypodiscretization}
The analysis of the equation
$
Tx=y,
$
with $T:X\to Y$ a~bounded linear operator
between some Banach spaces $X$ and~$Y$, consists in choosing
finite-dimensional  vector spaces
$X_N$ and~$Y_N$ and the corresponding embeddings $\imath_N$ and~$\jmath_N$:

 $$
    \bfig
\square/>`<-`<-`>/[X`Y`X_N`Y_N;T`\imath_N`\jmath_N`T_N]
 \efig
 $$
In this event, the equation
$$
 T_N x_N=y_N
$$
is viewed as a finite-dimensional approximation
to the original problem.

\section{\large\bf\color{blue}Hyperdiscretization}
\subsection{}
Nonstandard models yield the method of hyperapproximation

 $$\bfig
\square[E`F`E^{\scriptscriptstyle\#}`F^{\scriptscriptstyle\#};T`\varphi_{E}`\varphi_{F}`T^{\scriptscriptstyle\#}]
 \efig
 $$

Here $E$ and $F$ are normed spaces over the same
scalars, while
$T$ is a bounded linear operator from $E$ to~$F$, and
${}^{\scriptscriptstyle\#}$ symbolizes a nonstandard  hull.

\section{\large\bf\color{blue}The Hull of a Space}
Let ${}^\ast$ is the symbol of the
Robinsonian standardization. Let $(E,\|\cdot\|)$~be an internal normed space over
~${}^\ast\mathbb F$, with $\mathbb F:=\mathbb R;\mathbb C$.
 As usual, ~$x\in E$
is a~{\it limited\/} element provided that
$\|x\|$~is a limited real
(whose modulus  has a standard upper bound by definition).
If $\|x\|$ is an infinitesimal then  $x$
is also referred to as an {\it infinitesimal}.
Denote by $\ltd(E)$  and  $\mu (E)$  the external sets of limited elements
and infinitesimals of~$E$. The set $\mu (E)$ is the~ {\it monad\/}
of the origin in~$E$.
Clearly, $\ltd(E)$~is an external vector space over~$\mathbb F$, and
 $\mu(E)$~is a subspace of $\ltd(E)$.  Put
 $E^{\scriptscriptstyle\#} =\ltd(E)/\mu (E)$ and endow
 $E^{\scriptscriptstyle\#}$ with the
natural norm
 $  \|\varphi x\|:=\|x^{\scriptscriptstyle\#}\|:=\st (\|x\|)\in\mathbb F$ for all $x\in\ltd(E)$ Here
 $\varphi:=\varphi_E:=(\cdot)^{\scriptscriptstyle\#}:\ltd(E)\to E^{\scriptscriptstyle\#}$
 is the canonical homomorphism, and $\st$ takes the standard part of a limited real.
This $(E^{\scriptscriptstyle\#}, {\|\cdot\|})$ is an
external normed space called the
 {\it nonstandard hull\/} of~$E$.

\section{\large\bf\color{blue}The Hull of an Operator}
 Suppose now that
 $E$ and $F$ are internal normed spaces and
 $T:E \to F$ is an internal bounded linear operator.
The set of reals
 $ c(T):= \{C\in {}^\ast{\mathbb R}:\,(\forall x\in E)
 \| T x \| \leq C\|x \| \}
 $
 is internal and bounded. Recall that
 $\|T\|:= \inf c(T)$.
 If the norm
 $\|T \|$ of~$T$ is limited then the classical normative inequality
 $\| T x \| \leq \| T \|\,\|x\|$ valid for all $x \in E$,
 implies that
 $T({\ltd(E)})\subset \ltd(F)$
and
 $T({\mu (E)})\subset \mu (F)$.
 %$^{\scriptscriptstyle\#}$
 Hence, we may soundly define the
descent of~$T$ to the factor space $E^{\scriptscriptstyle\#}$
as the external operator
 $T^{\scriptscriptstyle\#}:E^{\scriptscriptstyle\#}\to
 F^{\scriptscriptstyle\#}$,
 acting by the rule
 $$
  T^{\scriptscriptstyle\#} \varphi_E x :=\varphi_F T x\quad (x \in E).
 $$
 The operator
 $T^{\scriptscriptstyle\#}$
 is linear  (with respect to the members of~${\mathbb F}$)
 and bounded; moreover, $\|T^{\scriptscriptstyle\#}\| =\st(\|T\|)$.
The operator
 $T^{\scriptscriptstyle\#}$ is called the
 {\it nonstandard hull\/}
of~$T$.

\section{\large\bf\color{blue}One Puzzling Definition}

Approximation of arbitrary function spaces
and operators by their analogs in finite dimensions, which
is {\it discretization},
matches the marvelous  universal
understanding of computational mathematics as the science of
finite approximations to general (not necessarily metrizable)
compacta. This revolutionary and challenging definition  was given in
the joint talk~ submitted by S.~L.~Sobolev, L.~A.~Lyuster\-nik,
and L.~V. Kantorovich at the Third All-Union Mathematical Congress
in~1956.

Infinitesimal methods  suggest a background, providing new
schemes for hyperapproximation of general compact spa\-ces.
As an approximation  to a compact space we may take
an arbitrary internal subset containing all standard elements
of the space under approximation.

\section{\large\bf\color{blue}State of the Art}

Adaptation of the  ideas of model theory
to optimization projects among the most important
directions of developing the synthetic methods of pure and applied
mathematics.  This approach yields new models of numbers,
spaces, and types of equations. The content expands of
all available theorems and algorithms. The whole methodology
of mathematical research is enriched and renewed, opening up
absolutely fantastic opportunities.
We can now use actual infinities and infinitesimals, transform
matrices into numbers, spaces into straight lines, and noncompact spaces into
compact spaces, yet having still uncharted vast territories of new knowledge.

\section{\large\bf\color{blue}Vistas of the Future}
Quite a long time had passed until the classical functional analysis
occupied its present position of the language of continuous mathematics.
Now the time has come of the new powerful technologies of model theory
in mathematical analysis. Not all theoretical  and applied mathematicians
have already gained the importance of modern tools and
learned how to use them. However, there is no backward traffic in
science, and the new methods
are doomed to reside in the realm of mathematics for ever and in a short
time they will become  as elementary and omnipresent in calculuses and
calculations as Banach spaces and linear operators.

\bibliographystyle{plain}

\end{document}